# Old and new results on normality

**Martine Queffélec**[1]

*Université Lille1*

**Abstract:** We present a partial survey on normal numbers, including Keane's contributions, and with recent developments in different directions.

## 1. Introduction

A central tool in the papers of Mike Keane is the ingenious and recurrent use of measure theory (probability and ergodic theory) in the various problems he investigated, paving the road for successors. So did he with normality, one of his favorite topics.

Let $q \geq 2$ be an integer. *A real number $\theta \in I := [0,1)$ is said to be normal to base $q$, if for every $m \geq 1$, every finite word $d_1 d_2 \ldots d_m \in \{0, 1, \ldots, q-1\}^m$ occurs in the $q$-adic expansion of $\theta$, $\sum_{n \geq 1} \epsilon_n(\theta) q^{-n}$, with frequency $q^{-m}$.* This means that every pattern appears infinitely often with the "good frequency" in the expansion of $x$.

A number $\theta \in I$ is said to be *simply normal to base $q$* if this property holds for $m = 1$. It is said to be *q (absolutely) normal* if it is normal to every base $q$, and it is *simply normal* if it is simply normal to every base $q$.

Normal numbers have been introduced in 1909 by E. Borel in [7], where he proved that–with respect to Lebesgue measure $m$–*almost every number is normal*. In fact, the digits $(\epsilon_j(\theta))_{j \geq 1}$ of $\theta$ in base $q$ behave as independent identically distributed random variables. Borel's result is thus a consequence of the Strong Law of Large Numbers, applied for each $q$ to disjoint blocks $\epsilon_{km+1}(\theta)\ \epsilon_{km+2}(\theta) \ldots \epsilon_{(k+1)m}(\theta)$, $k \geq 0$.

Normality to base $q$ can also be defined in terms of dynamics. The $q$-expansion of $\theta$ is obtained by iterating on $\theta$ the $q$-transformation $T_q$ defined on $I$ by $T_q(x) = qx \mod 1$; if $x_0 = x$ and $x_{n+1} = qx_n \mod 1$ for $n \geq 1$, then $\epsilon_{n+1}(x) = [qx_n]$. It is well known that the Lebesgue measure $m$ on $I$ is the unique absolutely continuous $T_q$-invariant measure, and that the dynamical system $(I, T_q, m)$ is ergodic. Borel's result is thus a consequence of the Birkhoff Ergodic Theorem.

So, normal numbers lie at the junction between two areas of expertise of Mike Keane, who gave illuminating proofs of both the Strong Law of Large Numbers and the Birkhoff Ergodic Theorem [19, 22], and wrote

<p align="center">"Measure-preserving transformations are beautiful"</p>

with the evident corollary

<p align="center">"Invariant measures are beautiful!"</p>

[1]Université Lille1, UMR 8546, France, e-mail: `martine@math.univ-lille1.fr`

## 2. Combinatorial point of view

From the definition, explicit $q$-normal numbers and non-normal numbers can be exhibited; for example Champernowne's number [11], whose expansion consists in concatenating all consecutive words to base $q$, is $q$-normal, for example

$$c_2 = 0.01000110110000010100111001011101111\ldots$$

when $q = 2$.

However, nobody knows whether classical arithmetical constants such as $\pi$, $e$, $\sqrt{2}$, $\zeta(n), \ldots$ are normal numbers to a fixed base, say $q = 2$. Even weaker assertions are unresolved. For example, it is not known whether there exist infinitely many 7s in the decimal expansion of $\sqrt{2}$; or whether arbitrarily long blocks of zeros appear in its binary expansion, i.e., $\liminf_{n\to\infty}\{2^n\sqrt{2}\} = 0$?

Returning to Champernowne's number $c_2$, the question arises whether $c_2$ is normal to any other base. Even though almost every number is normal to every base $q$, no explicit example is known of a number which is both $p$ and $q$-normal, with $p \neq q$. It is also not known how can one get the 3-expansion from the 2-expansion of some irrational number.

A plausible conjecture (suggested by Borel's work) is that *every irrational algebraic number is normal*. Already a partial answer to this conjecture such as *every irrational algebraic number is simply normal*, would be of great interest as we shall see.

An unsolved question from Mahler [26] is the following. *Suppose that $\alpha = \sum_1^\infty \frac{a_n}{2^n}$ and $\beta = \sum_1^\infty \frac{a_n}{3^n}$ are algebraic numbers, with $a_n \in \{0,1\}$. Are $\alpha$ and $\beta$ necessarily rational?*

If any irrational algebraic number would be simply normal, $\beta$ must be rational, and consequently the sequence $(a_n)$ is ultimately periodic; This would imply that $\alpha$ is rational as well, thus answering Mahler's question positively.

In recent investigations on normality, the statistical property of "good frequency" is replaced by the combinatorial property of "good complexity," asserting that $p(n)$, the number of words of length $n$ occurring in the $q$-expansion of $x$, is maximal, i.e., $p(n) = q^n$, $n \geq 1$, leading to the following conjecture.

**Conjecture.** *Every irrational algebraic number has maximal complexity function $p(n)$.*

Using a $p$-adic version of Roth's theorem, which appeared in previous work by Ridout [35], Ferenczi and Mauduit [12] obtained the following description of a class of transcendental numbers (in a revisited formulation by Adamczewski and Bugeaud).

**Theorem 2.1 (Ferenczi–Mauduit).** *Let $x$ be an irrational number, whose $q$-expansion begins with $0.U_n V_n^s$ for every $n \geq 1$, with*
  (i) $s > 2$;
  (ii) $|V_n|$, *the length of the word $V_n$, is increasing;*
  (iii) $|U_n|/|V_n|$ *is bounded.*
*Then $x$ is a transcendental number.*

As a corollary they get the transcendence of Sturmian numbers:

*If there exists $q$ such that the expansion of $x$ to base $q$ is a Sturmian sequence on $q$ letters, then $x$ is a transcendental number.*

Remembering that a Sturmian sequence on $q$ letters satisfies $p(n) = n + q - 1$ for $n \geq 1$, this can be viewed as the first result on complexity in this setting:

**Theorem 2.2.** *If $x$ is an irrational algebraic number, then its complexity function to base $q$ satisfies*

$$\liminf_{n \to \infty}(p(n) - n) = +\infty.$$

Roughly speaking Theorem 2.1 means that the sequence of digits begins with arbitrarily long prefixes of the form $W_n W_n W_n^\epsilon$; but these conditions eliminate some natural numbers; for example, the Thue-Morse sequence on $\{0, 1\}$ being without overlap of size $2 + \varepsilon$, and the associated Thue-Morse number as well. Involving a suitable version of the Schmidt Subspace Theorem, which may be considered as a multi-dimensional extension of previous quoted results by Roth and Ridout, Adamczewski and Bugeaud [1] improved this theorem by replacing in (i) $s > 2$ by $s > 1$. This proves as a by-result that the Thue-Morse number is transcendental. As a corollary they get

**Theorem 2.3 (Adamcewski–Bugeaud).** *If $x$ is an irrational algebraic number, then its complexity function to base $q$ satisfies*

$$\liminf_{n \to \infty} \frac{p(n)}{n} = +\infty.$$

Thus, every non-periodic sequence with sub-linear complexity gives rise to a transcendental number; this is the case with generalized Morse sequences [20, 21].

## 3. Metric point of view

We identify throughout $[0, 1)$ and the circle $\mathbf{T} = \mathbf{R}/\mathbf{Z}$, and let $M(\mathbf{T})$ be the set of bounded Borel measures on $\mathbf{T}$.

In this section we turn to the characterization of the sets $\mathcal{N}_q$, of $q$-normal numbers, and of $\mathcal{N} = \cap_q \mathcal{N}_q$ as well as their negligible complements, raising many connected questions.

1. Under what conditions on $p$ and $q$ (if any) does there exist a number normal to base $p$ but non-normal to base $q$?
2. In this case, what is the cardinality of $\mathcal{N}_p \cap \mathcal{N}_q^c$?
3. What is the Hausdorff dimension of such a set?

Another approach to normality will be useful. Recall that *a real sequence $(u_n)_n$ is said to be uniformly distributed mod 1 if*

$$\forall 0 \le a < b \le 1, \ \frac{1}{N}|\{n \le N \ ; \ \{u_n\} \in [a,b]\}| \to b - a,$$

*as $N \to \infty$, where $\{x\}$ is the fractional part of $x$. Equivalently, by the Weyl's criterion, if and only if*

$$\forall k \ne 0, \ \frac{1}{N} \sum_{n \le N} e_k(u_n) \to 0,$$

*where $e_k(x) := e^{2i\pi kx}$.* This means that the sequence of probability measures on $I$:

$$\frac{1}{N} \sum_{n \le N} \delta_{\{u_n\}} \to m \text{ weak}^*$$

where $\delta_x$ is the unit mass at $x$.

It was proved by Wall in [39] that $x \in I$ is normal to base $q$ if and only the sequence $(q^n x)_n$ is uniformly distributed mod 1. This point of view leads to the ergodic proof of Borel's theorem since, $m$-almost everywhere,

$$\frac{1}{N} \sum_{n<N} e_k \circ T_q^n \to \int_I e_k \, dm = 0.$$

### 3.1. Schmidt's results

In 1960, W. Schmidt solved questions 1. and 2. We write $p \sim q$ if there exist positive integers $r, s$ with $p^r = q^s$ (or $\frac{\log p}{\log q} \in \mathbf{Q}$), and $p \nsim q$ otherwise.

**Theorem 3.1 (W. Schmidt).** 1. *Assume $p \sim q$; then any number normal to base $p$ is normal to base $q$: $\mathcal{N}_p = \mathcal{N}_q$.*

*2. If $p \nsim q$, then the set of $p$-normal numbers which are not even simply normal to base $q$ has the power of the continuum.*

*Proof.* 1. The identity $\mathcal{N}_{p^r} = \mathcal{N}_p$ can be established by using the combinatorial definition. This gives the first assertion immediately.

An alternative proof of the inclusion $\mathcal{N}_{p^r} \subset \mathcal{N}_p$ goes as follows: if $x \in \mathcal{N}_{p^r}$ so are $px, p^2 x, \ldots, p^{r-1} x$, since $\sum_{n<N} e_k(p^{rn} p^j x) = \sum_{n<N} e_{kp^j}(p^{rn} x)$. The normality of $x$ to base $p$ is a consequence of the division algorithm.

For the reverse inclusion $\mathcal{N}_q \subset \mathcal{N}_{q^s}$ we can proceed as suggested in [3]. Suppose that $x \in \mathcal{N}_q$: $\nu_N = \frac{1}{N} \sum_{n \leq N-1} \delta_{\{q^n x\}} \to m$ weak*; we shall prove that

$$\nu_N^{(s)} = \frac{1}{N} \sum_{n \leq N-1} \delta_{\{q^{sn} x\}} \to m \text{ weak}^*$$

by proving that $m$ is the unique weak* limit point of the sequence $(\nu_N^{(s)})$. If $F \subset I$ is a closed set such that $m(F) = 0$, then by hypothesis $\nu_{sN}(F) \leq \varepsilon$ for $N \geq N_F$. The evident inequality $\frac{1}{sN} \sum_{n \leq N-1} \delta_{\{q^{sn} x\}}(F) \leq \nu_{sN}(F)$ yields that

$$\nu_N^{(s)}(F) \leq s\varepsilon, \ N \geq N_F.$$

Let $\sigma$ be a weak* limit point of $\nu_N^{(s)}$ in the weak* compact set of probability measures on $I$. From the above, $\sigma(F) = 0$ whenever $F$ is closed with $m(F) = 0$. For every Borel set $A$ with $m(A) = 0$, one has

$$\sigma(A) = \sup_{\substack{K \subset A \\ K \text{ compact}}} \sigma(K) = 0$$

since $m(K) = 0$, and $\sigma \ll m$. But $\sigma$ is $T_{q^s}$-invariant and $m$ is the unique absolutely continuous $T_{q^s}$-invariant measure. Hence $\sigma = m$ and $x \in \mathcal{N}_{q^s}$.

2. The proof of the second assertion, which is rather complicated in [37], has been simplified by many authors, first of them being Keane and Pierce [23]. The basic idea is that any $q$-invariant, ergodic and non absolutely continuous probability measure $\mu$ supports $\mathcal{N}_q^c$. By choosing $k$ with $\hat{\mu}(k) \neq 0$, and using the ergodic theorem, one gets

$$\frac{1}{N} \sum_{n<N} e_k(q^n x) \to \hat{\mu}(k) \neq 0 \ \mu-a.e.$$

and $\mu$-almost $x$ is in $\mathcal{N}_q^c$.

Keane and Pierce constructed a Cantor measure (depending on $q$ but not on $p$) supported by $\mathcal{N}_p \cap \mathcal{N}_q^c$. Later, generalizations have been obtained by Brown, Moran & Pierce [9] by employing Riesz products, by Feldman & Smorodinski [13] then Bisbas [6] with Bernoulli convolution measures, by Kamae [17] with specific singular measures; Host [14] improved the result in view of Fürstenberg's conjecture. Pollington [31] showed that the set $\mathcal{N}_p \cap \mathcal{N}_q^c$ has full Hausdorff dimension. We give here a revisited proof by Kamae [17] where $g$-measures, introduced by Keane [21], play an important role. □

Let us recall the notion of $g$-measures. Let $T$ be an $m$-to-1 covering transformation on a compact metric space $X$, and $g$ a positive continuous function on $X$ such that $\sum_{Ty=x} g(y) = 1$. We define the transfer operator $L$ on $C(X)$ by

$$Lf(x) = \sum_{Ty=x} f(y)g(y).$$

*A $g$-measure is a $T$-invariant Borel probability measure, which is an eigenmeasure to the maximal eigenvalue of the adjoint operator $L^*$.* Uniqueness holds whenever $g$ is strictly positive with Lipschitz regularity, and in this case $\mu$ is strongly mixing [21]. Natural examples appear in computing the correlation measure of finitely-valued $q$-multiplicative sequences of modulus 1 (which are nothing else than generalized Morse sequences).

Fix $q \geq 2$; consider $S_q(n)$ the sum of the $q$-digits of $n$ and $\rho := \rho_q$ the correlation measure of the sequence $\alpha(n) = e^{2i\pi\alpha S_q(n)}$, $\alpha \in \mathbf{R}$; $\rho$ is the strongly mixing $g$-measure $\prod P(q^n x)$, with $P(x) = \frac{1}{q}|\sum_{0 \leq k < q} \alpha(k)e(kx)|^2$ and $g = \frac{1}{q}P$. It is proved in [32] that *$\rho$ has the constant modulus property*; we need only a weaker assertion:

**Lemma 1.** *Denote by $\Gamma$ the dual group of $\mathbf{T}$ when acting multiplicatively. If $\chi_\rho \in \overline{\Gamma}_\rho$, the set of limit points for the topology $\sigma(L^\infty(\rho), L^1(\rho))$ of sequences $(\gamma_n) \subset \Gamma$, then $|\chi_\rho| = C$, $C$ constant.*

As explained above, $\rho(\mathcal{N}_q) = 0$. Now, for every $k \neq 0$,

$$\int \Big|\frac{1}{N}\sum_{n<N} e_k(p^n x)\Big|^2 d\rho(x) = \frac{1}{N^2}\sum_{\substack{n<N\\m<N}} \hat{\rho}(kp^n - kp^m)$$

and the property $\rho(\mathcal{N}_p) = 1$ for $p \not\sim q$ follows readily from the next one:

**Lemma 2.** *For every integer $a \neq 0$, $\lim_{\substack{n,m\to\infty\\n>m}} \hat{\rho}(a(p^n - p^m)) = 0$.*

*Proof.* Fix $a > 0$ and denote by $\chi_\rho \in \overline{\Gamma}_\rho$ a limit point in the topology $\sigma(L^\infty(\rho), L^1(\rho))$ of the sequence $(\gamma_{ap^n}) \subset \Gamma$, where $\gamma_k = e_k$. By Lemma 1, $|\chi_\rho| = C$. Since $|\chi_\rho|^2$ is a limit point of the sequence $(\gamma_{ap^n - ap^m})_{n>m}$, it suffices to prove that $C = 0$.

If $C \neq 0$, there exists $c \neq 0$ and $\gamma = \gamma_b$ with $|\int \chi_\rho \gamma_b \, d\rho| = c$; now, by extracting a subsequence if necessary,

$$\hat{\rho}(ap^n + b) \to \int \chi_\rho \gamma_b \, d\rho.$$

Let $0 < \varepsilon \ll 1$, and $\rho$ being continuous, let $m > 0$ be such that $|\hat{\rho}(m)| < \varepsilon c$. If

$q^{N-1} \leq m < q^N$, from hypotheses on $p, q$, there exists $n_k, m_k$ going to $\infty$ such that

$$\frac{p^{n_k}}{q^{m_k}} \to \frac{m}{aq^N} \in [0, 1[.$$

Thus, $ap^{n_k} = mq^{m_k-N} + r_k$ where $r_k = o(q^{m_k})$, and for $k$ big enough,

$$ap^{n_k} + b = mq^{m_k-N} + r_k + b, \ m_k > N, \ r_k + b < q^{m_k-N}.$$

By the strong mixing property,

$$\hat{\rho}(ap^{n_k} + b) \sim \hat{\rho}(m)\hat{\rho}(r_k + b).$$

But $\hat{\rho}(ap^{n_k} + b) \to c$ and $|\hat{\rho}(m)\hat{\rho}(r_k + b)| < \varepsilon c$, which yields a contradiction if such a $c \neq 0$ exists. This implies that $\chi_\rho = 0$ and the lemma is proved. □

This lemma has been proved in a different way by Kamae with help of a deep result from Strauss and Senge.

### 3.2. Sums of normal numbers

Observe first that a Borel set $E \subset \mathbf{T}$ with $m(E) > 1/2$ satisfies $E + E = \mathbf{T}$. This follows from the fact that for every $x \in \mathbf{T}$, $m(E \cap (x - E)) > 0$, hence there exist $u, v \in E$ such that $x = u + v$. Now if $m(E) = 0$ but $E$ supports a probability measure $\mu$ satisfying $\mu * \mu = m$, then every number in $\mathbf{T}$ can be expressed as the sum of **four** numbers in $E$ since $m(E + E) = 1$. But this is far from optimal.

It is well-known that $K_3$, the triadic Cantor set, is such that $K_3 + K_3 = \mathbf{T}$. However, if $\mu_3$ is the Cantor measure (the probability measure uniformly distributed over those points of $\mathbf{T}$ whose ternary expansions contain only 0's and 2's), $\mu_3 \notin M_0(\mathbf{T}) = \{\mu \in M(\mathbf{T}); \lim_{n \to \infty} \hat{\mu}(n) = 0\}$ and $\mu_3 * \mu_3$ is not even absolutely continuous with respect to the Lebesgue measure $m$ on $\mathbf{T}$.

More striking: Erdös gave an example showing that there exist arbitrarily thin sets, in the sense of Hausdorff measures, whose sum set contains an interval.

Now, if $E$ and $F$ have positive Cantor measure, $E + F$ has positive Lebesgue measure, since we have the following inequality

$$m(E + F) \geq \mu_3(E)^\alpha \mu_3(F)^\alpha \tag{3.1}$$

where $\alpha = \log 3 / \log 4$; this follows from the combinatorial one: If $A, B \subset \{0, 1\}^n$, then $|A+B| \geq |A|^\alpha |B|^\alpha$. (Note that $2\alpha = \dim m / \dim \mu_3$). The extension of (3.1) to the Cantor measure $\mu_m$ on $K_m$, the $m$-adic Cantor set, $m \geq 3$, has been established in [8]:

$$m(E_1 + \cdots + E_m) \geq \prod_{j=1}^{m} \mu_{m+1}(E_j)^{\alpha(m)}$$

where $\alpha(m) = \log(m+1)/m \log 2$.

Similar observations can be made about the negligible sets $\mathcal{N}(p) = \cap_{q \not\sim p} \mathcal{N}_q \cap \mathcal{N}_p^c$. Keane and Pierce have shown in [23] that $\mu_p(\mathcal{N}(p)) = 1$. Combining this with all previous remarks yields that every number in $\mathbf{T}$ can be expressed as the sum of $2(p-1)$ numbers in $\mathcal{N}(p)$. But various measures living on $\mathcal{N}(p)$ (as we saw) lead to different decomposition results.

By using Hausdorff measures, Pollington proved in [31] the following result.

**Theorem 3.2.** *For each $p \geq 2$, $\mathcal{N}(p) + \mathcal{N}(p) = \mathbf{T}$.*

*3.3. More about non-normal numbers*

So what more can be said about measures living on the negligible sets $\mathcal{N}_q^c$ ? Probability measures $\mu$ on **T** whose Fourier transform vanishes at infinity rather fast are still annihilating the sets $\mathcal{N}_q^c$. More precisely, in [25] the following was shown. *If $|\hat{\mu}(n)| \leq \phi(|n|)$ whenever $\phi$ is a non-increasing function such that $\sum \frac{\phi(n)}{n \log n} < \infty$, then $\mu(\mathcal{N}_q^c) = 0$.* This is a consequence of the following Strong Law of Large Numbers for weakly correlated random variables, applied to $X_n = e_{kq^n}$.

**Theorem 3.3 (Davenport–Erdös–LeVeque).** *Let $\mu$ be a probability measure on $E$, $(X_n)$ a sequence of bounded random variables and $S_N = \frac{1}{N} \sum_{n=1}^{N} X_n$. If*

$$\sum_{N \geq 1} \frac{1}{N} \int |S_N|^2 d\mu < \infty,$$

*then $S_N \to 0$ $\mu$-a.e.*

On the other hand, Lyons constructed a probability measure concentrated on $\mathcal{N}_2^c$ with $|\hat{\mu}(n)| \leq \phi(|n|)$, where $\phi$ is a non-increasing function such that $\sum \frac{\phi(n)}{n \log n} = \infty$. It is a convolution type measure with respect to which the binary digits as random variables are independent by blocks. His sharpest result is the following theorem.

**Theorem 3.4 (Lyons).** *There exists a probability measure $\mu \in M_0(\mathbf{T})$ such that*

$$\limsup_{K \to \infty} \frac{1}{K} \sum_{1}^{K} \mathbf{1}_{[0,1/2]}(2^k x) = 1 \quad \mu - a.e.$$

Then arised the following question: Does there exist $x$ and $I$, closed arc of length $1/2$, such that $\{2^n x\} \in I$ for all $n \geq 0$?

## 4. Normality and continued fractions

Every real number admits a regular continued fraction expansion, which is more handy than the $q$-expansions, and furnishes the best rational approximations. Although it is difficult to get one expansion from the other, there exist classes of examples for which both are known, but no example of normal number with bounded partial quotients has been exhibited yet.

We recall some facts concerning the regular continued fraction algorithm (RCF).

Given an irrational number $x \in I$, $x$ is the limit of the infinite sequence of rational numbers

$$\frac{P_k(x)}{Q_k(x)} = \cfrac{1}{a_1(x) + \cfrac{1}{a_2(x) + \cdots \cfrac{1}{a_{k-1}(x) + \cfrac{1}{a_k(x)}}}},$$

where the integers $a_k(x) \geq 1$ if $k \geq 1$; we write shortly $\frac{P_k}{Q_k} = [0; a_1, \ldots, a_k]$, $x = [0; a_1, a_2, \ldots]$ where the *partial quotients* of $x$, $a_1, a_2, \ldots$, are uniquely defined.

The shift on the RCF expansion is conjugated to an expanding transformation, the Gauss map, defined on $X = [0,1] \backslash \mathbf{Q}$ by $Tx = \dfrac{1}{x}$ mod 1. The partial quotients

$a_n = a_n(x)$, $n \geq 1$, are then given by

$$a_1(x) = [\frac{1}{x}], \ a_{k+1}(x) = a_k(Tx).$$

The probability measure $\mu$ on $X$ (called the *Gauss measure*) defined by $\mu(A) = \frac{1}{\log 2} \int_A \frac{dx}{1+x}$ for every Borel-set $A$, is preserved by $T$ and is $T$-ergodic. Furthermore, $\mu$ is the unique absolutely continuous $T$-invariant probability measure.

Applying Birkhoff's ergodic theorem, we get that

$$\lim_{n \to \infty} \frac{1}{n}(a_1(x) + \cdots + a_n(x)) = +\infty \quad \mu - a.e.$$

Since $\mu$ is equivalent to Lebesgue measure, the *set BAD of numbers in $[0,1)$ with bounded partial quotients has zero Lebesgue measure*.

A significant difference between adic- and RCF expansions is the behavior of the digits as random variables; the $a_n(x)$ are no more independent identically distributed variables on $(X, \mu)$, but quasi-independent in the following sense.

**Theorem 4.1 (Philipp).** *Consider the $\sigma$-algebras $\mathcal{A}_k = \sigma(a_1, \ldots, a_k)$ and $\mathcal{B}_k = \sigma(a_{k+1}, \ldots)$. Then there exists $\theta \in ]0, 1[$ such that*

$$|\mu(A \cap B) - \mu(A)\mu(B)| \leq \theta^n \mu(A)\mu(B)$$

*for every $A \in \mathcal{A}_k$ and $B \in \mathcal{B}_{k+n}$, $k \geq 1$, $n \geq 1$.*

For complementary notations and results on the regular continued fraction expansion we refer to [16] and [5].

In [24], a class of normal numbers with explicit RCF expansion is exhibited; this combines two ideas.

1. *Every number of the form $\sum_k \frac{1}{A_k}$, with $A_k$ integers $\geq 2$ such that $A_k^2$ divides $A_{k+1}$, has an explicit RCF expansion.*

This is a consequence of classical identities around symmetry, for example [28, 38], if $\frac{P}{Q} = [0; a_1, \ldots, a_N]$, $a_N \geq 2$, $x \geq 1$,

$$\frac{P}{Q} + \frac{(-1)^N}{xQ^2} = [0; \underbrace{a_1, \ldots, a_{N-1}}, a_N, x-1, 1, a_N - 1, \underbrace{a_{N-1}, \ldots, a_1}].$$

In this way, Mahler proved that the transcendental number $\sum_{k=0}^{\infty} \frac{1}{u^{2^k}}$, $u \geq 2$, is in BAD.

2. *Every number $\alpha$ of the form $\sum_k p^{-\lambda_k} q^{-\mu_k}$, where $(\lambda_k), (\mu_k)$ are increasing sequences of integers $\geq 1$, $p, q \geq 2$ relatively prime integers and $\mu_k \geq p^{\lambda_k}$, is normal to base $q$.*

This is due to sharp estimates on the distribution of digits in periodic fractions and a consequence of the structure in base $q$ of the numbers $\alpha$. It can be seen that the $q$-expansion of an $\alpha$ looks like $0^a W_0^{s_0} W_0' W_1^{s_1} W_1' \ldots W_n^{s_n} W_n' \ldots$, $W_n'$ being a prefix of the word $W_n$, $|W_n| = t_n$ the order of $q$ mod $p^{\lambda_n}$ and $|W_n^{s_n} W_n'| = \mu_{n+1} - \mu_n$.

**Theorem 4.2 (Korobov).** *Numbers of the form $\alpha = \sum_k p^{-2^k} q^{-p^{2^k}}$ with $p, q$ relatively prime integers, are normal to base $q$ with explicit RCF expansion.*

Of course in these cases $x_k = A_{k+1}/A_k^2 \to \infty$ and such numbers $\alpha \notin$ BAD.

Montgomery [29] in his book posed the following problem: Does there exist normal numbers with bounded partial quotients ?

In 1980, Kaufman [18], using the structure of the sets $F(N) = \{x \in [0,1); x = [0; a_1, a_2, \ldots]$ with $a_i \leq N \ \forall i \geq 1\}$, obtained the following deep result.

**Theorem 4.3 (Kaufman).** *For $N \geq 3$ $F(N)$ carries a probability measure $\mu$ such that $|\hat{\mu}(t)| \leq c|t|^{-\eta}$ for a certain $\eta > 0$ (Kaufman's measure).*

When Kaufman's paper appeared, R.C. Baker observed in [29] that combining the above result with a corollary of the Davenport–Erdös–Leveque's theorem (already used by Lyons for non-normal sets), yields the existence of normal numbers in BAD.

More generally, if $\mathcal{A}$ is a finite alphabet of integers $\geq 1$, $|\mathcal{A}| \geq 2$, and

$$F(\mathcal{A}) = \{x \in [0,1);\ x = [0; a_1, a_2, \ldots] \text{ with } a_i \in \mathcal{A}\ \forall i \geq 1\},$$

In [34] we improved Kaufman's result as follows.

**Theorem 4.4.** *If $\dim_H(F(\mathcal{A})) > 1/2$, $F(\mathcal{A})$ supports a Kaufman's measure and contains infinitely many normal numbers. In particular there exist infinitely many normal numbers with partial quotients $\in \{1, 2\}$.*

Note that no explicit normal numbers in BAD have been constructed yet.

Now, it is natural to ask about the existence of RCF-normal numbers. Such a number $x$ is defined by the property that *for every $m \geq 1$, every finite word $d_1 d_2 \ldots d_m \in \mathbf{N}^{*m}$ occurs in the RCF expansion of $x$, with frequency $\mu([d_1 d_2 \ldots d_m])$, the Gauss measure of the cylinder $[d_1 d_2 \ldots d_m]$*. Once more from ergodicity (or quasi-independence), almost all $x \in [0,1)$ are normal for the continued fraction transformation and once more constructing one RCF-normal number raises difficulties. A successful proceedure has been carried out by Adler, Keane and Smorodinsky [2].

**Theorem 4.5.** *Let $Q_n$ be the ordered set of all rationals in $[0,1)$ with denominator $n$. The RCF expansion obtained by concatenating the RCF expansions of the numbers in $Q_2, Q_3, \ldots$ leads to an RCF-normal number.*

Hence the following problems: How to construct a number normal with respect to both adic- and Gauss transformations? Does there exist an RCF-normal number with low complexity relative to an adic-expansion?

## 5. More investigation

### 5.1. Topological point of view

The sets $\mathcal{N}_q^c$, $\mathcal{N}^c$ are small from the metric point of view but big from the topological one. This has been known for a long time, for example, a result from Helson and Kahane [15] goes as follows. *For each $q \geq 2$, $\mathcal{N}_q^c$ intersects every open interval in a uncountable set.* Also, if $q \geq 2$, $0 \leq r < q$, $x \in [0,1)$ and $N_n(r, x)$ is the number of occurrences of $r$ in the first $n$ terms of the $q$-expansion of $x$, then the set of limit points of the sequence $(N_n(r, x)/n)$ is $[0, 1]$, for all $x$ but a set of first Baire category; as a consequence *the set $\mathcal{N}$ is in first Baire category* [36].

Real numbers with non-dense orbit under $T_q$ ("$q$-orbit") are in fact very interesting. If $I$ is an open arc of **T** of length $> 1/2$, the set of $x \in \mathbf{T}$ whose 2-orbit does not intersect $I$ must be finite, but there exist $s \in [0, 1/2]$ whose 2-orbit is infinite and whole contained into $[s, s+1/2]$; this furnishes an answer to the question raised at the end of section 4; in fact we have a complete description of these numbers $s$ (see [4, 33]).

**Theorem 5.1 (Bullett–Sentenac).** *Let $s = \sum_{j \geq 1} \frac{\varepsilon_j}{2^j}$ be an irrational number in $[0, 1/2]$. Then following assertions are equivalent.*

*1. The closed 2-orbit of $s$ lies into $[s, s + 1/2]$.*

*2. The sequence of digits $(\varepsilon_j)_j$ is a characteristic sturmian sequence on the alphabet $\{0, 1\}$.*

### 5.2. Generalizations

Generalization of normality to non-integer bases or to endomorphisms leads to very interesting questions (see [3]). But this would be another story ...

### Acknowledgments

Thanks to the organizers of this colloquium in honor of Mike, thanks to the referee, and my gratitude to Mike himself for transmitting his enthousiasm through various lectures he gave us many years ago.